\PassOptionsToPackage{hyphens}{url}
\documentclass{notices}
\usepackage{amsfonts,amssymb,amsmath,amscd}
\usepackage{xcolor}
\usepackage[USenglish]{babel}
\usepackage{csquotes}
\usepackage{graphicx}
\usepackage{enumitem}
\usepackage{framed}
\usepackage{hyperref}

\usepackage{tikz}
\usetikzlibrary{positioning}
\usetikzlibrary{arrows}
\usepackage{tablefootnote}

\usepackage[normalem]{ulem}

\let\oldhref\href
\renewcommand{\href}[2]{\oldhref{#1}{\texttt{#2}}}

\usepackage{abstract}

\title{Making mathematical online resources FAIR:\\
at the example of small phylogenetic trees}

\author{
  Tabea Bacher%
  \affil{TB works at the Max Planck Institute for Mathematics in the Sciences. Her email address is 
  \url{bacher@mis.mpg.de}.}
  \and
  Marina Garrote-L\'opez%
  \affil{MGL is a researcher at the Universitat Pompeu Fabra in Barcelona. Her email address is 
  \url{marina.garrote@upf.edu}.}
  \and
  Christiane G\"orgen%
  \affil{CG (corresponding author) is a mathematical data steward at Leipzig University. Their email address is 
  \url{goergen@math.uni-leipzig.de}.}
  \and
  Marius Neubert%
  \affil{MN is a PhD student at Leipzig University. His email address is 
  \url{marius.neubert@math.uni-leipzig.de}.}
}

\begin{document}

\maketitle 



\section{Introduction}

\begin{figure*}[t]
  \centering
  \includegraphics[width=0.295\textwidth]{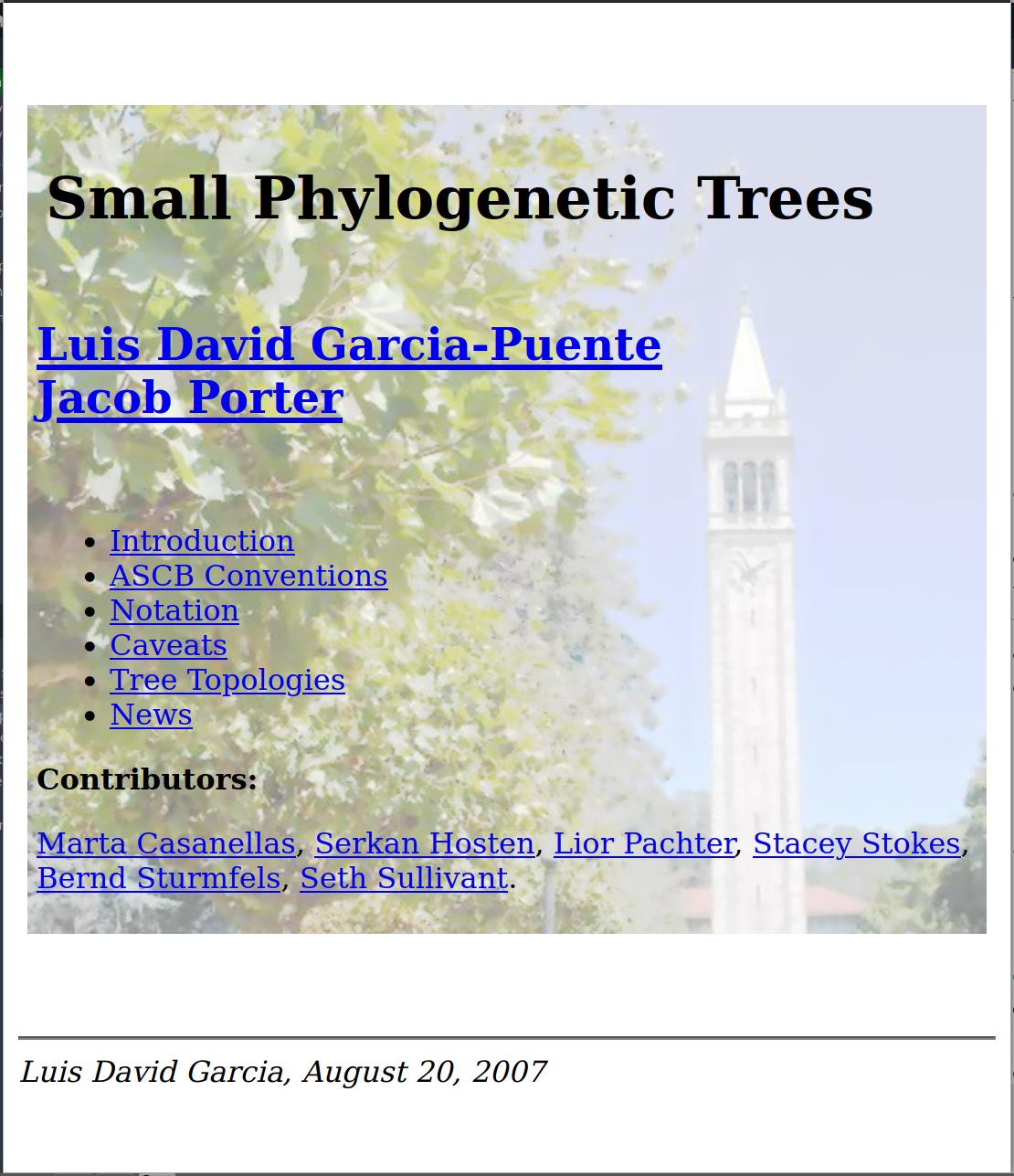}
  \quad
  \includegraphics[width=0.295\textwidth]{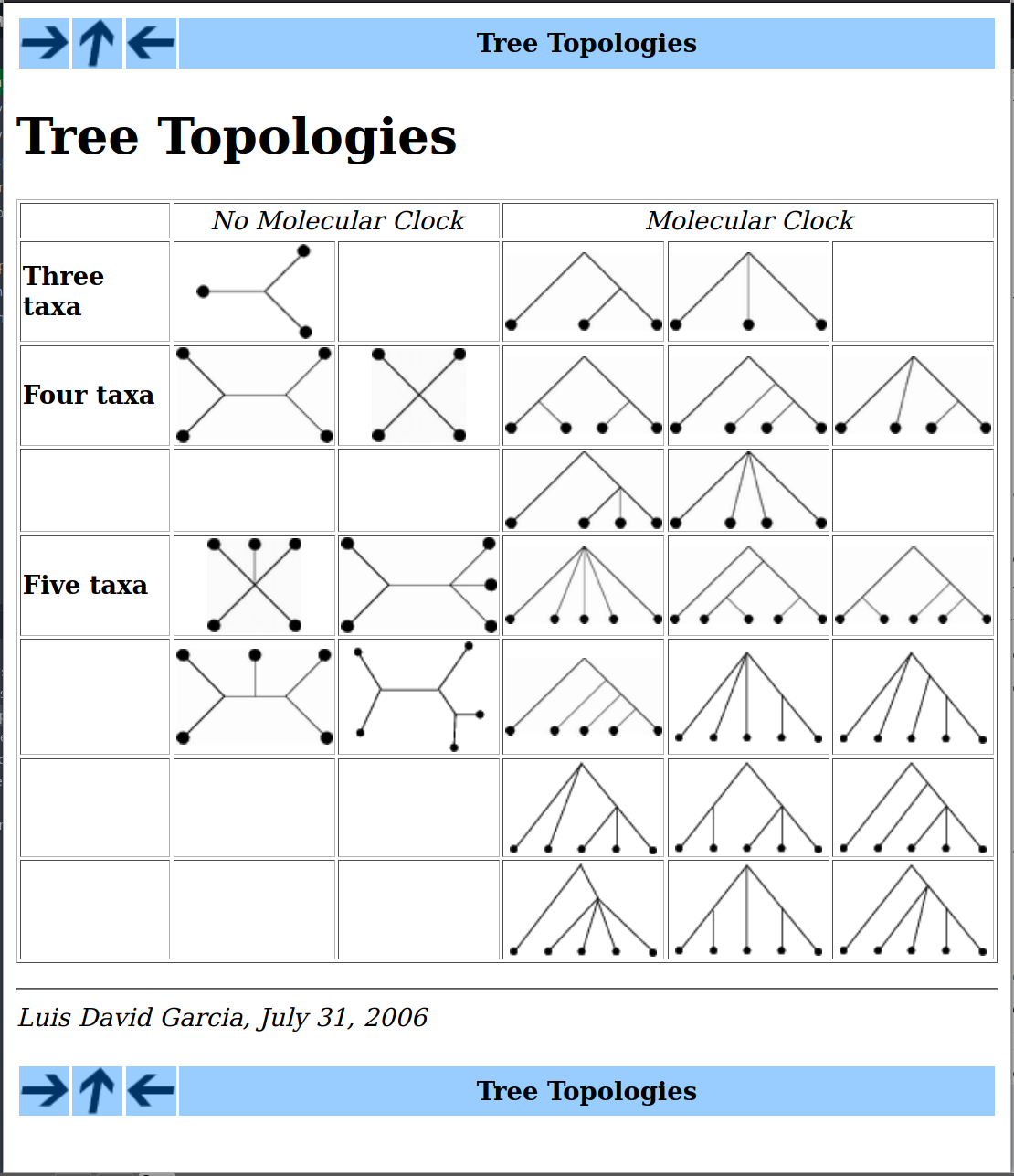}
  \quad
  \includegraphics[width=0.295\textwidth]{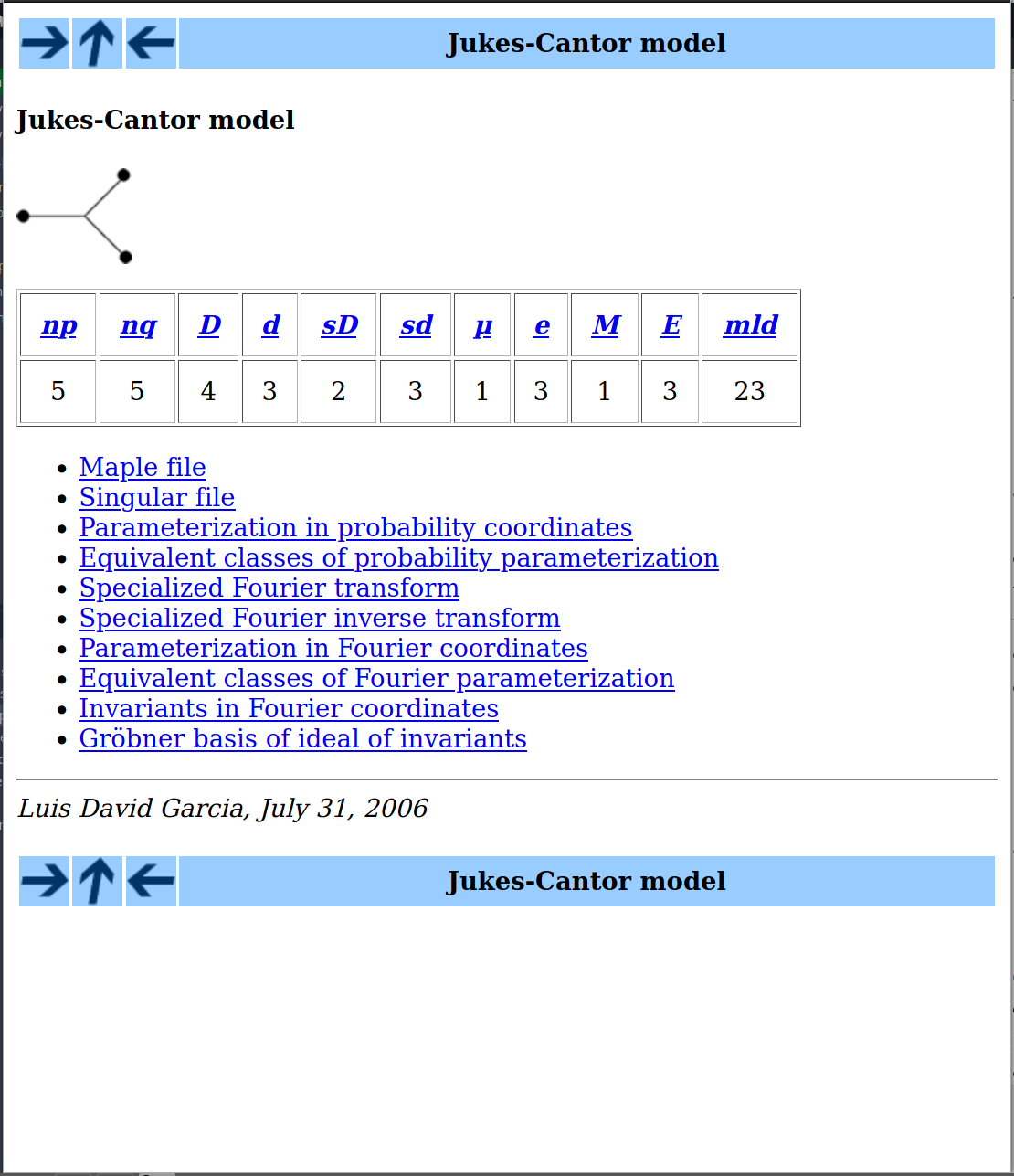}
  \caption{Screenshots of Small Phylogenetic Trees of January 2025. The first (a) shows the landing page, the second (b) an overview of all trees considered, and the third (c; after having navigated through two more subpages) data on the three-leaf phylogenetic tree with Jukes-Cantor model.}\label{fig:spt}
\end{figure*}

The mathematical library \emph{Small Phylogenetic Trees} \cites{zbMATH06811710,SmallPhylogeneticTrees}, first compiled in the early 2000s, is a significant digital resource in the field of \emph{algebraic phylogenetics}. It catalogs tables of algebraic invariants of evolutionary models on phylogenetic trees up to five species. 
These invariants provide insight into the algebro-geometric structure of these statistical models and may be used for phylogenetic inference. We give a brief introduction to their mathematical description in Section~\ref{sec:algPhylo}.
Small Phylogenetic Trees supplements the invariants with code, raw computational output files, and a documentation of the terminology used, see~Fig.~\ref{fig:spt}. 
Twenty years later, it remains an actively used online resource for researchers in algebraic phylogenetics but no longer fully meets contemporary digital standards.

This paper analyzes its current standing and presents a generalizable approach detailing how to preserve the library for the future. Rather than introducing new mathematical results, it serves as a methodological report and case study drawing on expertise from academic research-data managers and mathematicians, offering practical guidance for others seeking to revitalize and sustain similar projects.

\medskip

In the last two decades, digital technology has transformed mathematical research.
In the early 2000s, the digital landscape for mathematics grew with the rapidly-developing public internet. 
%
%
The now-famous \emph{On-Line Encyclopedia of Integer Sequences} (OEIS), available at \href{https://oeis.org}{oeis.org}, had been online for about a decade after two decades as a book, while the equally-famous LMFDB (\emph{L-functions and modular forms database}) at \href{https://lmfdb.org}{lmfdb.org}, {also based on earlier books,} was still just a vague idea.
The arXiv preprint server was yet to reach half a million articles (cf.\,\href{https://arxiv.org/stats/main}{arxiv.org/stats/main}), and prominent mathematicians like Gowers and Tao had just taken up their blogs, \href{https://gowers.wordpress.com}{gowers.wordpress.com} and \href{https://terrytao.wordpress.com}{terrytao.wordpress.com}. The version-control system \texttt{git} had only been recently launched. 
In this {environment, the goal of} Small Phylogenetic Trees to \enquote{make available in a unified [digital] format} the known {computational} results in algebraic phylogenetics was modern and visionary.

Today, the situation has changed significantly. Digital Object Identifiers (DOIs), introduced at the turn of the century, have become standard for citing publications, in many aspects taking over the role of ISBNs for identifying books. {Git repository services} like \texttt{GitHub} and institutional \texttt{gitlab} are now indispensable for collaborative coding. 
Web technologies now allow to dynamically present a lot of data on a single page. Powerful open-source software for symbolic computation rivals commercial alternatives. 
Knowledge graphs can link mathematical results directly to their formalizations and related papers, and mathematical notation may be displayed and annotated using {\LaTeX} and, e.g., MathML.

\medskip

Digital preservation today is guided by the FAIR Principles \cite{Wilkinson.etal.2016}, which advocate for \emph{Findability, Accessibility, Interoperability}, and \emph{Reusability} of digital research results across the sciences; see Section \ref{sec:FAIR} for more details. These principles do not demand for the use of the newest and shiniest web tool for presentation, replacing old standards by ever new ones---much more sensibly, they guide our decisions on selecting methodologies that preserve the long-term value and usability of digitized results.

By today's standards, the Small Phylogenetic Trees library appears dated. Is it FAIR? Are its results verifiable? 
We will analyze these questions in Section \ref{sec:analysisSPT}.
In this case study, we treat Small Phylogenetic Trees as a prototypical, valuable mathematical resource, rich with significant and timeless results---and in need of modernization to ensure its survival. 
Our approach is threefold:
\begin{enumerate}[nosep]
    \item We launch a new, straightforward presentation of the original website content; built as a collaborative project and findable at \href{https://algebraicphylogenetics.org}{algebraicphylogenetics.org} \cite{AlgebraicPhylogenetics}.
    \item We select a single software environment capable of computationally verifying every result on the website, and we use it to develop a comprehensive package for algebraic phylogenetics, \texttt{AlgebraicPhylogenetics} \cite{AlgebraicPhylogeneticsJL}.
    \item We compile this text as the third component: tracing our methodology, decision-making processes, and implementation details. This crucial aspect of the project was not documented in the original book chapter on Small Phylogenetic Trees, but we consider it vital to enable the community to follow by example.
\end{enumerate}

Our primary objectives are: firstly, to enhance the library's adherence to the FAIR Principles and, secondly, to ensure the long-term reproducibility of its mathematical results. Thus, in Section \ref{sec:plan} we document  the details of our strategic approach, and in Section \ref{sec:achivements} the outcomes achieved. {Finally,} in Section \ref{sec:learnings} we derive guidelines to support the sustainable modernization of other, similar digital mathematical resources.

\section{Background}

Although theoretical proofs will always be at the core of mathematics, modern research increasingly relies on computers to perform complex calculations and to analyze large sets of information. 
Online mathematical libraries such as the OEIS, LMFDB and Small Phylogenetic Trees have become the modern successors to look-up tables and classifications in book form. They are easier to update than their traditionally printed counterparts, they consolidate information in one place rather than spread out over a series of articles, and they serve as active tools where researchers can test conjectures against a wide range of cases, compute examples, and access hard-to-find results, thus saving significant time, effort, and energy. 
%

\subsection{Algebraic phylogenetic models}\label{sec:algPhylo}

Small Phylogenetic Trees self-describes as \enquote{mainly a repository of algebraic information of different evolutionary models based on trees with a small number of leaves.} Specifically, it focuses on nucleotide substitution models describing the evolution of DNA sequences, {composed of the nucleotides \texttt{A}, \texttt{C}, \texttt{G}, and \texttt{T},} on phylogenetic trees of five or fewer species, or taxa. Figure~\ref{fig:spt}b shows a screenshot of such trees. 
%
%
These models assume that DNA sequences evolve according to a Markov process on a given phylogenetic tree $T$. A random variable representing a nucleotide is associated with each node of the tree, and the process is governed by two key parameter sets: (1) a probability distribution at the root of $T$, and (2) a $4\times 4$ transition matrix associated with each edge.
This setup defines a polynomial map $\phi_T$ that sends the collection of free parameters, namely the root distribution and the transition matrices, to the joint probability distribution $p = (p_{\tt AA\ldots A}, p_{\tt AA\ldots C}, \ldots, p_{\tt TT\ldots T})$ of observing specific nucleotides at the leaves of $T$. {We refer to \cite{sullivant2018algebraic} for a more comprehensive introduction.}

{These models can be viewed as algebraic statistical models, since the set of all possible distributions {$p\in\text{Im }\phi_T$} that can arise from a phylogenetic tree $T$ lies within an algebraic variety $V_T$, the Zariski closure of the image of $\phi_T$.}
This variety is defined by a set of polynomials that vanish for any valid distribution generated by the model; they belong to the ideal $I(V_T)$, and are known as \emph{phylogenetic invariants}. 

While the theoretical properties of phylogenetic-tree models have been extensively studied \cites{Sturmfels2005, AllmanRhodes2008}, 
computing  {phylogenetic, as well as} \emph{algebraic invariants} such as the dimension or degree of $V_T$ or its singular locus, is a significant challenge. These computations often rely on having a Gröbner basis for $I(V_T)$ at hand, notoriously computationally intensive even for trees with relatively few taxa. This computational bottleneck is precisely the problem that the catalog of Small Phylogenetic Trees addresses: the user does not have to perform extensive computations or browse papers, but may simply look up algebraic invariants already computed and made available by trustworthy peers. 

The library focuses on a particular class of nucleotide substitution models known as \emph{group-based} models, such as the \emph{Cavender-Farris-Neyman}%
, \emph{Jukes Cantor} and \emph{Kimura} models. 
These allow for a Fourier transformation to simplify the parametrization of $V_T$. After this linear change of coordinates, the models correspond to toric varieties, making their phylogenetic invariants significantly easier to describe and compute \cites{HendyPenny1989, Sturmfels2005}.
{The website catalogs a range of these algebraic invariants, see~Fig.~\ref{fig:spt}c; for each model and tree, it contains data like the dimension (D) and degree (d) of $V_T$ and its singular locus (sD and sd, respectively) or the dimension of the linear space containing the probability distributions (np). }
Alongside the algebraic invariants, the website lists the parametrization $\phi_T$ in both standard and Fourier coordinates, and the defining phylogenetic invariants.
%
For instance, we see that the maximum likelihood degree, a measure of parameter estimation complexity, for a 3-leaf tree under the Jukes-Cantor model is a manageable 23 while the corresponding value for the 5-star tree remains unknown Without diving into computational details, a researcher can learn from Small Phylogenetic trees that the corresponding toric variety for the 5-star tree has dimension 6 and degree 115, lies in a 27-dimensional linear space and is defined by 175 polynomial equations of degree 4 and 5. 
%

\subsection{FAIR, reproducible, and open mathematics}\label{sec:FAIR}

The FAIR Principles are well established and widely accepted for assessing the sustainability of digital research artifacts such as papers, code, data or software. They comprise fifteen criteria \cite[Box~2]{Wilkinson.etal.2016} that help both human and machines determine
\begin{itemize}[nosep]
    \item[--]  whether a resource and its description is \emph{findable} online, and persistently so,
    \item[--] who is allowed to, and technically how the resource and its description can be \emph{accessed},
    \item[--] whether the resource is \emph{interoperable} with its environment,
    \item[--] whether, as well as under which licensing, it may be \emph{reused}.
\end{itemize}
FAIRness shall not be confused with reproducibility: while the former may be read as \emph{available and usable}, we read the latter as \emph{verifiable}. Thus, FAIR mathematics is not necessarily good or understandable mathematics. But a FAIR research output is a necessary first step towards assessing its reproducibility. Likewise, the accessibility \emph{A} in FAIR should not be confused with \emph{open access}. Although there is a strong open-access culture in mathematics (with preprints freely available under CC-BY licensing on \href{https://arxiv.org}{arxiv.org}), FAIR accessibility is fulfilled if data and metadata are machine readable, for instance via an API. The data itself does not need to be open, as long as its metadata is available.

Whilst the FAIR Principles are the key guiding principles for data management across the sciences, they are currently not widely known or consciously applied within the mathematics community. There are of course flagship projects of mathematical databases like the OEIS and LMFDB whose online presentation lives up to the same high standards the community applies to their theoretical findings---but we were unable to find general advice on how to make Small Phylogenetic Trees as FAIR (or reproducible) as these big players.

The National Research Data Initiative \href{https://nfdi.de}{nfdi.de} in Germany is likely the largest national effort to promote FAIR research data, with its mathematics-specific consortium MaRDI, cf.\, \href{https://mardi4nfdi.de}{mardi4nfdi.de}, providing insight on how to apply the FAIR Principles specifically to mathematical artifacts. The NFDI is part of international initiatives such as the European Open Science Cloud and the Research Data Alliance.
Current guidelines from general scientific funding bodies primarily focus on the FAIRness of typical research data types like datasets, publications, and software, offering little to no guidance specifically tailored to online mathematical libraries like Small Phylogenetic Trees. The MaRDI white paper \cite{MaRDI.2023} and the analysis \cite{Conrad.etal.2024} by MaRDI authors are, as far as we are aware, the only current resources providing some guidance and initial assessment of the FAIRness of some mathematical resources; but no suggestions on how to FAIR online mathematical libraries either. This paper aims to close that gap.

\section{FAIRifying small trees}
\subsection{Analyzing the website}\label{sec:analysisSPT}

Careful thought and effort has gone into Small Phylogenetic Trees. Before assessing its current state of FAIRness and reproducibility, we provide an overview of the library's design and content.

The website is structured like a traditional catalog, reflecting early-2000s web practice. Its landing page sports a table of contents: mathematical background, notation and tree-labeling conventions, errata of the print version, and a comprehensive table of small trees considered. It links to news items dated 2005 up to 2007, suggesting active maintenance in that period. 
The table of trees then enables a user to dive deeper: navigating through a subpage with all nucleotide substitution models considered for a given tree, to a subpage with algebraic invariants and code for each particular choice of tree and model. 

While the library's structure is timeless in its logic, its full potential is hindered by this nested design. Together with the use of unintuitive subpage identifiers (such as numbers rather than descriptive names) navigation is difficult for new users and broken author links evoke an outdated first impression.

Now what can we say about these data?

\begin{description}
    \item[Findability]%
%
Small Phylogenetic Trees is hosted by the main maintainer's academic institution and can be found through a general search engine, even if the direct link in the original paper \cite{zbMATH06811710} is  broken due to a change in affiliation and the database title does not match precisely. Researchers may use both the (hardly memorable) URL to refer to Small Phylogenetic Trees or cite the corresponding book chapter.
\item[Accessibility]
The website uses the \texttt{http} standard protocol. There is no mechanism in place to download data or code other than manually navigating through subpages, using copy and paste. Outdated file formats, such as \texttt{gif}, do not allow users to access the input of computations, such as the tree $T$, in an automatized way.
\item[Interoperability] 
%
The code is supplied in two programming languages which are standard in the field, commercial \texttt{Maple} and open-source \texttt{Singular}, enabling the user to make a choice compatible with their own preference and budget.
Interoperability is hindered
by the absence of input/output specifications, detailed software requirements, and versioning, and by poor documentation. How variable names in the code relate to the underlying theoretical concepts is in places hard to trace.
\item[Reusability]
%
Each subpage of the resource includes a timestamp (typically 2006) and an author (typically key maintainer Luis David Garcia Puente), providing some versioning and accountability.

There is no license attached to the data or code, making it legally not usable. 
However, all contributing authors are linked to their personal homepages, cf.\, Fig.~\ref{fig:spt}a, and may be contacted for reuse. 
\end{description}

We may assume that the lack of licensing is due to a typical pitfall for newcomers in research-data management, namely assuming no license would imply free usage. This is not the case but understandable given mathematics' culture of openly sharing results. 

The interoperability in FAIR is also interlinked with how reproducible the data turn out to be.

\begin{description}
    \item[Reproducibility]
    Standard algebraic-phylogenetics definitions, illustrations, and naming conventions help a user familiar with the field to trace a mathematical argument.
    However, important data such as the polynomial parametrization $\phi_T$ of $V_T$ and its phylogenetic invariants are hard-coded. That is, {the provenance of} key mathematical results {is} undocumented. This absence of methodological and computational transparency hinders direct verification.
\end{description}

Especially for newcomers, considerable effort is needed to make results on Small Phylogenetic Trees interoperable with and verifyable on their own systems. 
So at present, Small Phylogenetic Trees fulfills some of the FAIR Principles and is partly reproducible with sufficient domain knowledge and coding skills. But does this mean it is not trustworthy? No.
The resource's credibility is established by its attention to authorship, the standing of authors in the community, underlying papers which have passed peer review, and a broad and educated user base in the community. What we can do to improve the current state of the library is to make every aspect---from display and data retrieval to reuse and reproducibility---significantly more straightforward.

\subsection{Making a plan}\label{sec:plan}

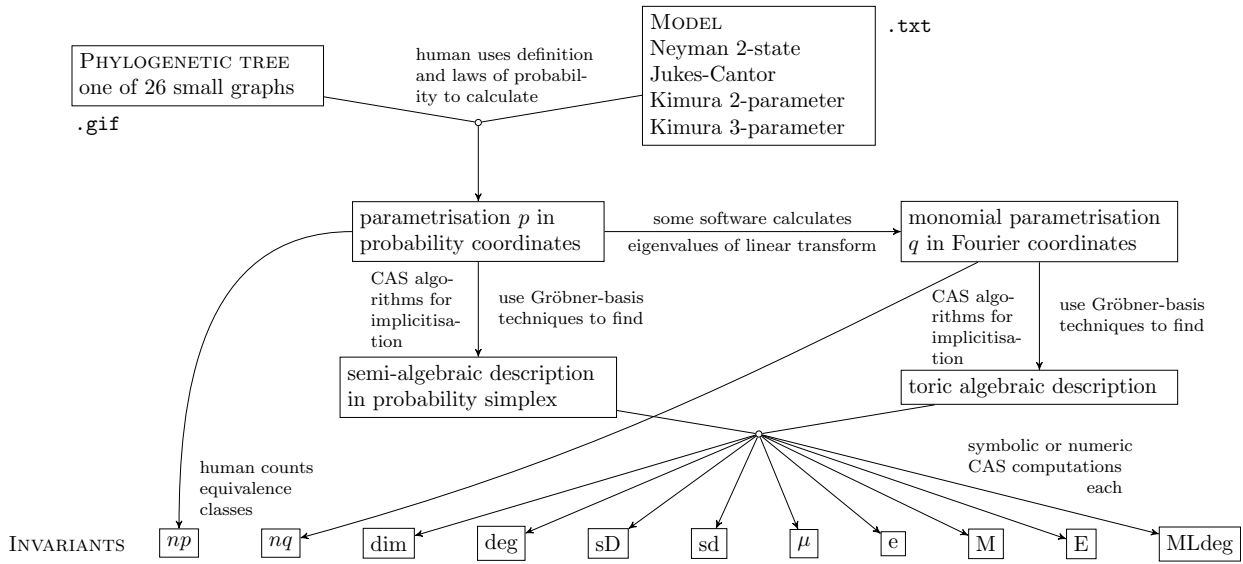
\begin{figure*}[t]
    \resizebox{1.02\textwidth}{!}{%

\begin{tikzpicture}[xscale=3,yscale=2.4,>=stealth']

\node[text width=38mm,draw,rectangle] (tree) at (1.2,0.2) {%
    \textsc{Phylogenetic tree $T$} 
    (one of 26 small graphs)};
\node[below left of=tree,xshift=-10mm,yshift=-0mm] {\footnotesize \texttt{.gif}};

\node[text width=6.9cm,draw,rectangle] (model) at (4.3,0.3) {%
    \textsc{Nucleotide Substitution Model}\\
    - {Neyman 2-state} {(Cavendar-Farris-Neyman)} \\    
    - Jukes-Cantor \\    
    - Kimura 2-parameter \\    
    - Kimura 3-parameter};
\node[below right of=model,xshift=12mm,yshift=-6mm] {\footnotesize heading in \texttt{.html}};

\node[text width=3.8cm,draw,rectangle] (probcord) at (0.7,-1) {%
    Parametrization $p$ in \\
    probability coordinates};
\node[below right of=probcord, xshift=10mm,yshift=-0mm] {\footnotesize \texttt{.txt}};

\node[text width=4.2cm,draw,rectangle] (fouriercord) at (5.3,-1) {%
    Monomial parametrization \\ $q$ in Fourier coordinates};
\node[below right of=fouriercord, xshift=12mm,yshift=-0mm] {\footnotesize \texttt{.txt}};

\node[text width=2.3cm] (phyloinv) at (-0.35,-2.2) {\textsc{Phylogenetic}\\ \textsc{Invariants}};
\node[text width=3.5cm,draw,rectangle, dashed] (algprob) at (0.85,-2.2) {%
    Algebraic description of $V_T$
    in probability coordinates $p$};

\node[text width=5.1cm,draw,rectangle] (algfourier) at (5.3,-2.2) {%
    Algebraic description of the toric $V_T$
    in Fourier coordinates $q$};
\node[below right of=algfourier, xshift=16mm,yshift=-0mm] {\footnotesize \texttt{.txt}};

\node[text width=2.3cm,draw,rectangle] (fouriertransf) at (3, -1) {%
   Fourier linear transformation};
\node[below right of=fouriertransf, xshift=2mm,yshift=-0mm] {\footnotesize \texttt{.txt}};

\node[text width=3cm,draw,rectangle] (invfouriertransf) at (3.4, -2.2) {%
   Inverse Fourier linear transformation};
\node[below right of=invfouriertransf, xshift=5mm,yshift=0mm] {\footnotesize \texttt{.txt}};

{
\node[text width=2.3cm] (prop) at (-0.35,-3.5) {\textsc{Algebraic}\\ \textsc{Invariants}};
\node[draw,rectangle] (np) at (0.4,-3.5) {$np$};
\node[draw,rectangle,right=of np] (nq) {$nq$};
\node[draw,rectangle,right=of nq] (dim) {dim};
\node[draw,rectangle,right=of dim] (deg) {deg};
\node[draw,rectangle,right=of deg] (sD) {sD};
\node[draw,rectangle,right=of sD] (sd) {sd};
\node[draw,rectangle,right=of sd] (mu) {$\mu$};
\node[draw,rectangle,right=of mu] (e) {e};
\node[draw,rectangle,right=of e] (M) {M};
\node[draw,rectangle,right=of M] (E) {E};
\node[draw,rectangle,right=of E] (MLdeg) {MLdeg};
}

\node[draw,circle,inner sep=1pt] (merge) at (2.7,-0.3) {};
\node [below left= of merge,text width=3cm, font=\footnotesize,yshift=11mm,xshift=-15mm] {human uses definition and laws of probability {to calculate}};
\node [below right= of merge,text width=5.5cm, font=\footnotesize,yshift=8mm,xshift=30mm] {{human uses} explicit parametrization (see \hspace*{6mm} \cite[Sec 15.3]{sullivant2018algebraic})  {to calculate}};

\draw [-] 	(tree) -- (merge);
\draw [-] 	(model) -- (merge);
\draw [->] 	(merge) edge node [below, sloped, pos=0.3, font=\footnotesize] {{(\emph{missing})}} (probcord); 
\draw [->] 	(merge) edge node [below, sloped, pos=0.3, font=\footnotesize] {{(\emph{missing})}} (fouriercord);


\draw [->] 	(probcord) -- node [above, font=\footnotesize,yshift=0mm] {CAS computation}
node [below, font=\footnotesize,yshift=-0mm] {(\emph{missing})} (fouriertransf);
\draw [->] 	(fouriertransf) -- node [above, font=\footnotesize,yshift=0mm] {CAS computation}
node [below, font=\footnotesize,yshift=-0mm] {(\emph{missing})} (fouriercord);

\draw [->, dashed] 	(probcord) -- (algprob);

\draw [->] 	(fouriercord) -- node [left,font=\footnotesize,xshift=1mm,text width=2.5cm,pos=0.6] {CAS computation\\ \hspace*{12mm}(\emph{missing})}
node [text width=2.5cm, right, font=\footnotesize,xshift=0.1cm, yshift=-2mm] {use Gr\"obner basis---based techniques to compute the implicitization} (algfourier);


\draw [->] 	(algfourier) -- (invfouriertransf);
\draw [->] 	(invfouriertransf) -- node [above,font=\footnotesize,xshift=7mm,pos=0.6] {Output of computations}
node [below, font=\footnotesize,xshift=3mm] {in \texttt{Singular} and \texttt{Maple}} (algprob);

\draw [<->, dashed] (fouriertransf) -- (invfouriertransf);

\draw (probcord) edge[out=230, in=110, ->,font=\footnotesize,swap,text width=2cm] (np);
\draw[->] (fouriercord)  ..controls (0.3,-2) and (2.3,-2).. (nq);

\node[font=\footnotesize,text width=2cm,above=of nq,xshift=-5mm,yshift=-9mm] {human/software counts equivalence classes};

\node[draw,circle,inner sep=1pt] (common) at (3.2,-2.8) {};

\draw [->, dashed] 	(algprob) -- (common);
\draw [-] 	(algfourier) edge node [below, sloped, pos=0.5, font=\footnotesize] {{(\emph{missing})}} (common);

\draw[->] (common) -- (dim);
\draw[->] (common) -- (deg);
\draw[->] (common) -- (sD);
\draw[->] (common) -- (sd);
\draw[->] (common) -- (mu);
\draw[->] (common) -- (e);
\draw[->] (common) -- (M);
\draw[->] (common) -- (E);
\draw[->] (common) -- node [text width=3.5cm, right,font=\footnotesize,yshift=0mm,xshift=20mm] {symbolic or numeric\\ \hspace*{5mm} CAS computations \hspace*{18mm} each} (MLdeg);

\end{tikzpicture}
%
    }
  \caption{This is a data layout plan, showing which mathematical and computational objects appear on Small Phylogenetic Trees, detailing their interdependencies and provenance. We use it to decide how to build the software package and which information to display on the website. CAS is an abbreviation for computer algebra system.}\label{fig:datalayout}
\end{figure*}

We decide to address concerns raised in our analysis above carefully and separately by splitting the computation of algebraic invariants from their display: by supplementing a new website with a new software package. This enables us to both modernize the presentation and, independently, to provide proof in a single software package that stores all necessary information (definitions, functions, inputs, and outputs) in one place. The software we decide on is the computer-algebra system \texttt{OSCAR} \cite{OSCAR}, based on \texttt{Julia}. It can save mathematical objects together with their context: e.g., for phylogenetic invariant, a specification of the ring in which the polynomials live. The software is further able to call both symbolic and numerical methods, so there is no need to switch systems when trying out different approaches to compute {algebraic} invariants.

Usability, both technically as well as visually, is our foremost concern; and our website-and-software approach allows us to differentiate between levels of user need and skill. The website itself shall be easy to navigate for the interested pure mathematician or newcomer with no coding skills. For those in need of code and detail, each result displayed shall be linked back to a paragraph in the software package's documentation as well as to the respective lines of underlying code and a serialization file; ideally supplemented also with citations of theory. All technical detail is thus made available for researchers to embed  computational results directly into their own work. This interconnectedness of data (invariants, formulae, formalizations, code, and publications) enables both human users and machines to grasp mathematical results and generate knowledge.

For efficient teamwork, to clarify priorities, and to make sure we address both FAIRness and reproducibility, we draft three plans: one data layout plan and one research-data management plan both for the software and the website. 
The first, displayed in Fig.~\ref{fig:datalayout}, illustrates what data and data relations we find in the {Small Phylogenetic Trees} library, and how these build on each other when tracing (reproducing) results. 
The latter plans, in Fig.~\ref{fig:rdmps}, build on this. They formulate a data-handling strategy for the software and the website, clarifying how to treat all occurring files and data formats during and after the project, ensuring they are handled FAIR.
\medskip

\begin{figure*}
\begin{framed}%
\emph{Software research data management plan:}
\indent We reuse data in the form of articles and code, of type \texttt{.pdf} and \texttt{.txt}. We generate new \texttt{Julia} code, to be integrated into the research-software project \texttt{OSCAR}, for findability and interoperability, as well as for high-quality maintenance and easy access. We serialize our computational output in \texttt{JSON}-based \texttt{.mrdi} files \cite{MR4786726}, and verify it by comparison to existing results via a \texttt{Perl} script. We produce \texttt{.tex} and \texttt{.pdf} for an article presenting the software package. Newly generated data does not exceed 1GB in size.

\indent We use the automated \texttt{Julia} feature for code documentation, following \texttt{OSCAR}'s established style guide.
{\texttt{AlgebraicPhylogenetics} is framed inside the new \texttt{OSCAR} module in \href{https://docs.oscar-system.org/dev/Experimental/AlgebraicStatistics/introduction/}{\texttt{AlgebraicStatistics}}, developed in parallel and in collaboration with Antony Della Vecchia, Tobias Boege and Ben Hollering.} 
For version control, \texttt{AlgebraicPhylogenetics} and all supplementary files will be on Ben Hollering's GitHub instance\footnotemark\newcounter{f1}\setcounter{f1}{\value{footnote}}\ 
with new versions merged into the \texttt{OSCAR} system's GitHub\footnotemark\newcounter{f2}\setcounter{f2}{\value{footnote}}. BH's repository can be edited by all authors of this project. Local copies of all branches provide backup during the project runtime. The accompanying article will be collaboratively written using the platform \href{https://overleaf.com}{overleaf.com}.

\indent Our software will be pulled into \texttt{OSCAR}, increasing longevity and functionality as compared to a standalone algebraic-phylogenetics project; available under the GNU Public License, Version 3.0+.

\indent During project runtime, MN sets up a first version of the code and implements serialization, TB is responsible for testing, CG for documentation, and MGL oversees the project, manages communication with \texttt{AlgebraicStatistics.jl}, and takes the lead in code development. After the project is completed, the \texttt{OSCAR} group continues maintenance, and MGL will keep the main responsibility for further development. We apply for funding for MN's position for six months. For open-source software and in-kind human resources, no any additional funds are needed.
\end{framed}

\begin{framed}%
\emph{Website research data management plan:}
We reuse files building the website \href{https://fanography.info}{fanography.info}, of type \texttt{.py}, \texttt{.html}, \texttt{.css}, \texttt{.md}, and \texttt{.txt}. Since they are not supplied with a license, we ask the author (Pieter Belmans) for approval. We generate new data of the same type, supplemented with illustrations in \texttt{.png} and \texttt{.tex}, embeddable using \texttt{tikzjax}, and formulae in \texttt{MathJax} for automated findability and interoperability. We newly generate \texttt{.yaml} files from \texttt{AlgebraicPhylogenetics}, and use \texttt{Flask} in \texttt{Python}, to display these data. All files use the same naming convention, \verb|<datatype>_<numerical_tree_id>-<alphabetic_model_id>.<filetype>|, for interoperability and reusability. New data will not exceed 1GB.

The website is collaboratively developed on CG's Github\footnotemark\newcounter{f3}\setcounter{f3}{\value{footnote}}\ for versioning and storage, local copies serve as backup. The code is supplied with commentation. All authors of this article have access to the repository. Slack and email are used for communication. For each release, we fork the project to \href{https://gitlab.mis.mpg.de/goergen/algebraicphylogenetics.org}{gitlab.mis.mpg.de/goergen/algebraicphylogeneticsorg} for website hosting for as long as one of the authors is affiliated to the institute. Both repositories are freely available under the MIT license. The domain \href{https://algebraicphylogenetics.org}{algebraicphylogenetics.org} is secured for ten years.

During active development, CG and MGL jointly take the lead, coordinating contributors and content, MN provides the \texttt{.mrdi}-interface between website and software, and TB supervizes design as well as licensing. After the first release, CG maintains the website and ensures compliance with FAIR Principles; MGL coordinates community development for both software and website. We allocate funds to hire a student assistant to build the first prototype. Hosting and domain acquisition is an in-kind contribution.
\end{framed}
\vspace{-0.4cm}
\caption{Focusing on the mathematics first and on their display second, we make two research data management plans: one for the software project and one for the corresponding website. }\label{fig:rdmps}
\end{figure*}
\footnotetext[\value{f1}]{\href{https://github.com/bkholler/OSCAR.jl}{github.com/bkholler/OSCAR.jl}}
\footnotetext[\value{f2}]{\href{https://github.com/oscar-system/Oscar.jl}{github.com/oscar-system/Oscar.jl}}
\footnotetext[\value{f3}]{\href{https://github.com/chgoergen/algebraicphylogeneticsorg}{github.com/chgoergen/algebraicphylogeneticsorg}}

We start with the \emph{data layout plan} in Fig.~\ref{fig:datalayout} to better understand the structure of Small Phylogenetic Trees. For each phylogenetic tree and each phylogenetic model listed at the very top of the chart, we can compute a number of algebraic invariants, abbreviated and listed at the very bottom. The workflow from input to output requires a human to infer modeling assumptions from a picture, and write those into some computer algebra system used to compute phylogenetic invariants as well as the rest of the algebraic invariants. 
The overview of all the steps involved in this process enables us to identify which ones need to be formalized in a comprehensive software package: we mark these \enquote{missing} in the plan.

The chart serves as a basis for making the following decisions for improvements:%
\begin{itemize}[nosep]
    \item[--] to display input (model and tree) and output (algebraic invariants) on the new website in one single big table,
    \item[--] to display in-between computations, such as $\phi_T$ or the phylogenetic invariants, on an individual subpage for each tree/model combination, linked to directly from within the respective cell of the main table
    \item[--] for the software package to
    \begin{itemize}[nosep]
    \item[$\circ$] decide and implement the input structure 
    \item[$\circ$] write functions for every computation of relevance,
    \item[$\circ$] store the output in an interoperable format.
    \end{itemize}
\end{itemize}
\medskip

Having, thanks to the data layout plan, understood which data need be reproducible, we then address our goal of developing future-compatible, FAIR online mathematics. This is based on two pillars. First, we make \emph{machine readable} every file and every piece of mathematics we handle: to be actionable for semantic web technologies and thus ready for coming research-data infrastructure tailored to mathematics\footnote{Compare for instance the services developed by MaRDI.}. And second, we make \emph{human readable} every such piece of data via extensive documentation: thus preserving its value in case technology ceases to exist.

Now to achieve that goal, we plan which data types we generate and how these are treated along their life cycle. For instance, we prefer open-source data types like \texttt{.tex} for the report and \texttt{.mrdi} for computational results, and we prefer data files which are understandable and readable by humans, like \texttt{.yaml} for the data displayed on the website. The \texttt{.yaml} files are, at least regarding their mathematical content, duplicates of the \texttt{.mrdi} files. We need both because the latter are not yet printable. Both are (reproducibly) generated automatically from our code. Moreover, we choose \texttt{git} technology for versioning and collaborative coding; in particular, GitHub for active code development, and an academic gitlab instance for long-term website hosting.

The \emph{research data management plans} documenting these decisions are displayed in Fig.~\ref{fig:rdmps}.
They are set up after a phase of exploration and before we take up any of the actual work, and they follow an established structure: employed, e.g., by the data-management community and funding bodies\footnote{\href{https://https://new.nsf.gov/funding/data-management-plan}{new.nsf.gov/funding/data-management-plan}}, to address management as well as practical responsibilities. Rather than being a mere bureaucratic addition to the core mathematical work, they make it easy for us to focus because all possible distractions are clarified in advance: in particular, the order in which to perform tasks, whom to attribute which role to, how to keep track of progress---and, most importantly, how to handle and preserve results long term.


\subsection{Analyzing our achievements}\label{sec:achivements}

At the time of writing, an experimental version of our software is available \cite{AlgebraicPhylogeneticsJL} and so is version 0.1 of \href{https://algebraicphylogenetics.org}{algebraicphylogenetics.org}. 
\begin{figure*}[t]
  \centering
  \includegraphics[width=0.45\textwidth, trim=0 -228mm 0 0]{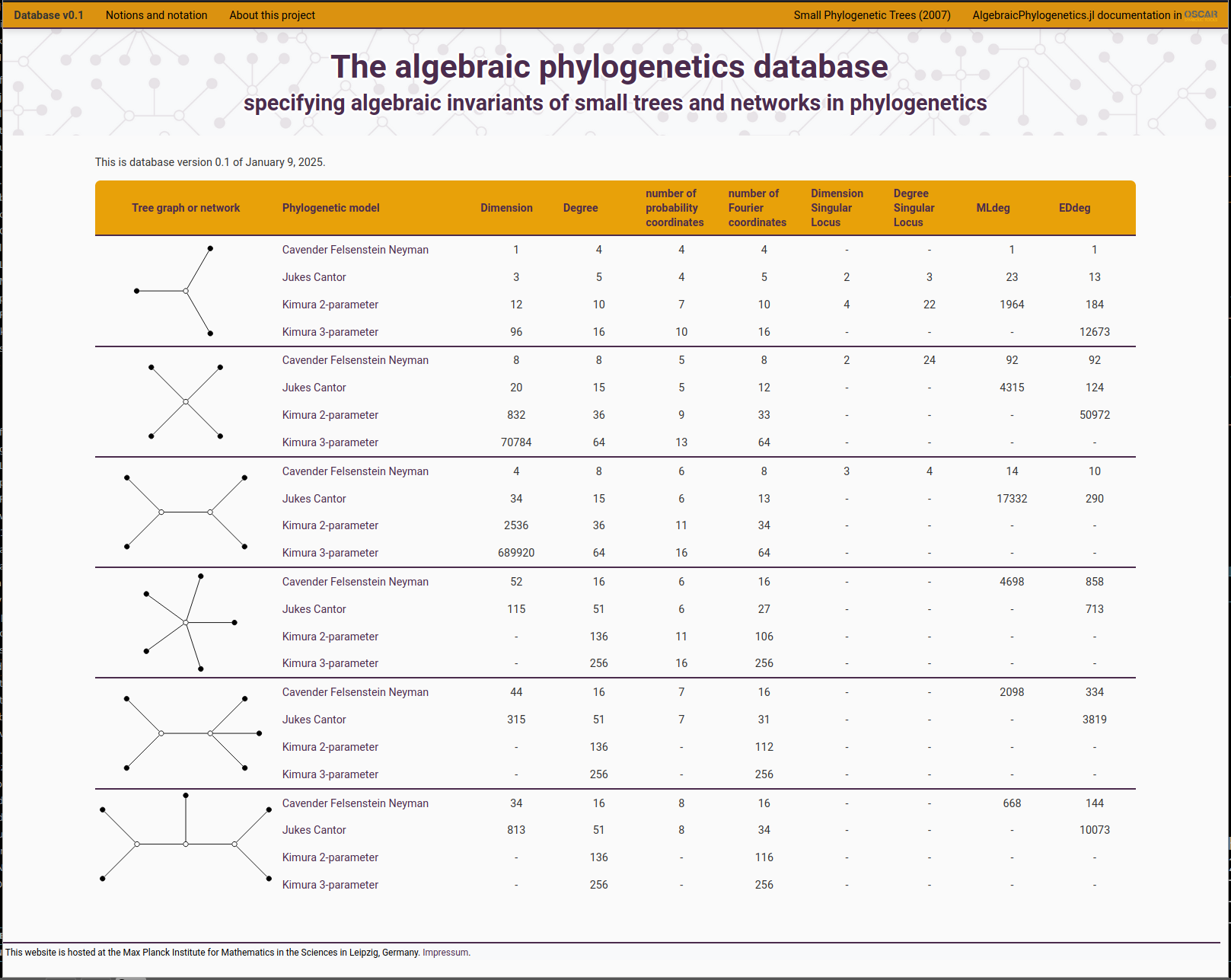}
  \qquad
  \includegraphics[width=0.45\textwidth]{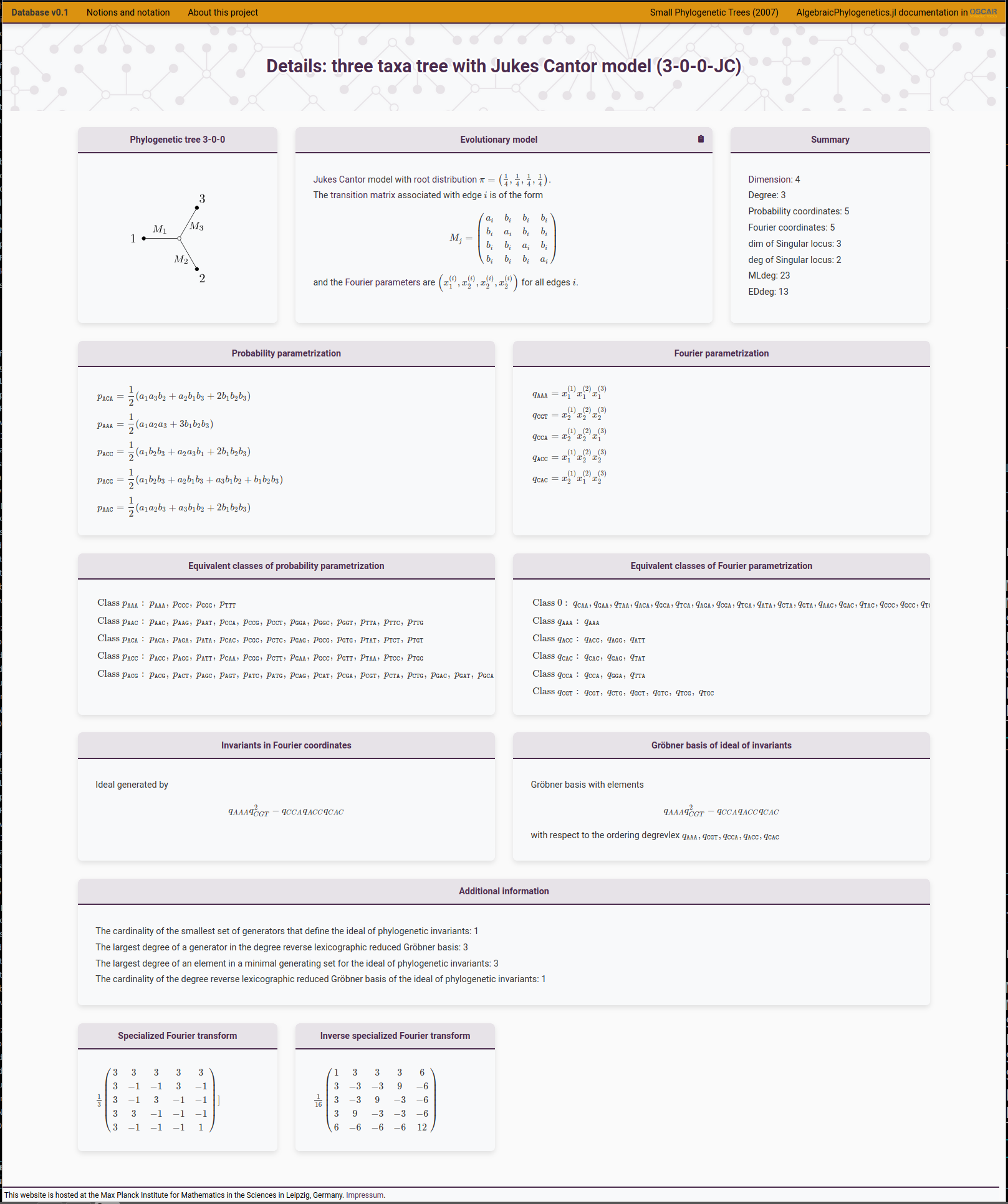}
  \caption{This is a screenshot of \href{https://algebraicphylogenetics.org}{algebraicphylogenetics.org} under development in January 2025: on the left is the landing page, on the right data on the three-leaf phylogenetic tree with Jukes-Cantor model.}\label{fig:apd}
\end{figure*}

In Fig.~\ref{fig:apd}, we present two screenshot of our new website. The landing page consists on a large table indexed by tree and nucleotide-substitution model, containing the algebraic invariants. This allows researchers in phylogenetics to quickly access knowledge.
At the top left of the page, there are two links directing users to subpages on notation and terms of use, respectively. Two additional links in the top right lead to Small Phylogenetic Trees and to {\texttt{AlgebraicPhylogenetics}}.
From the main table, users can directly navigate to the subpage of interest to find additional algebraic information about specific trees and models, and to integrate relevant computations and formulae into their own work.


We made improvements in the handling of mathematical data as well as in their metadata description. These are, in particular, of the following kind.

\begin{description}
    \item[Findability] The new domain and name, the \emph{Algebraic Phylogenetics Database} at \href{https://algebraicphylogenetics.org}{algebraicphylogenetics.org}, refer directly to the mathematical field, and we state the respective MSC identifiers 62R01, 65H10, 92D15, 68W30, and 14Q30 for findability. The homepage has its own citation guideline \cite{AlgebraicPhylogenetics} which is distinct from text reports (like this one) and from the software citation \cite{AlgebraicPhylogeneticsJL} for added precision and appropriate credit. For long-term preservation and citability, each past and current version of the website is stored at \href{https://zenodo.org}{zenodo.org}. as a snapshot of our GitHub repository. These means increase findability both for the human user and for automated online search engines.
    
    Additionally, we link all authors to their orcids and personal webpages, and we set up a maintainers email address. This enables users to directly and correctly identify lasting contact points. Relevant publications in algebraic phylogenetics are linked via DOIs for persistence. 
    Moreover, \texttt{AlgebraicPhylogenetics} benefits from OSCARs SWHID \cite{OSCAR} \href{https://softwareheritage.org}{software heritage.org} for long-term traceability. 
    \\

    \item[Accessibility] Metadata for each combination of phylogenetic tree and model are gathered in files which are both machine- and human-readable for immediate accessibility and long-term understandability. In particular, this information is stored in \texttt{.yaml} files which are hosted on GitHub and accessible via GitHub's REST API. Formulae inside these files are typeset in \LaTeX{} as is custom in mathematics, machine-readably rendered on the website using \texttt{MathJax}, and trees coded in \texttt{tikz} and displayed using \texttt{tikzjax}. Each mathematical object is additionally available in open-source \texttt{.mrdi} files based on \texttt{JSON} which can be loaded in \texttt{Julia}, thus containing valuable mathematical metadata again in machine-readable and human-readable format.  License statements on the GitHub projects underlying the website as well as the software give access to any interested user.%

    \item[Interoperability] 
    Regarding interoperability for the human, we reiterate standard notation, document our computational workflow and provide usage examples for the software package to facilitate integration into external projects.

    Maximum machine interoperability is achieved by choosing the open-source computer algebra system \texttt{OSCAR} which integrates \texttt{ANTIC}, \texttt{GAP}, \texttt{polymake}, and \texttt{Singular}; all of these widely used in computer algebra. Here, the \texttt{Julia} basis has the potential to facilitate the combination of both numeric and algebraic calculations in algebraic phylogenetics. The \texttt{.mrdi} format can provide both an interface between website and software as well as serialization: allowing to meaningfully save and load computations, much preferably to copying and pasting \LaTeX{} without context.%
    \item[Reusability]
    Extensive documentation both on the website itself and externally in \texttt{OSCAR} make it easy for the newcomer to dive into details. Permissive open licenses are chosen and clearly stated for displayed results, computer-algebra software, and the website build.
    
    Every result is further linked to a copyable {\LaTeX} display of formulae and to GitHub or \texttt{OSCAR} code performing relevant computations, making results traceable and boosting reproducibility and trustworthiness, as is good scientific practice.
\end{description}

In this way we have created an online resource that is as FAIR as our technical expertise, time constraints, and currently existing mathematics research-data infrastructure have enabled us to build.
\medskip


Regarding the reproducibility of results in the library, we have made the following improvements.

\begin{description}
    \item[Reproducibility] 
Coding in open-source software enables the users to check every function and data type we defined. The software \texttt{AlgebraicPhylogenetics} captures all computations in one place, as opposed to missing intermediate steps or results spread out over many different files. We provide extensive code documentation in \texttt{OSCAR}'s style guide, using the \texttt{Julia} documentation feature, and we link all data displayed back to the code. Data files are easily identifiable thanks to the new file-naming conventions. Any outsider can now simply redo computations and verify the data comparing input to output. Serialization provides extra reproducibility by seamlessly fitting into existing \texttt{OSCAR} environments. Versioning of our software enables practitioners to check compatibility. Our software passes the newly-formulated guidelines for writing and reviewing software in computer algebra \cite{Hanselman.2025}.
\end{description}

Both \texttt{Julia} and \texttt{OSCAR} are existing academic projects with big communities. This enables us to draw on developers' expertise, acquire new users, greatly improving chances of long-term maintenance as opposed to a standalone project, and crucially increase the chance that our work will continue to be FAIR and reproducible for the future.

We anticipate that the database will continue to grow with current interest in the phylogenetics community. In future versions, we plan to incorporate algebraic data related to phylogenetic networks, as well as additional nucleotide substitution models, such as algebraic time-reversible models. This opens the door to external contributions both through direct involvement in development and by sharing emerging interests, improvements, or updates that should be reflected in the database and accompanying software.


\section{Learnings and guidelines}\label{sec:learnings}
\subsection{Lessons learned}

Being a rather large team with heterogeneous expertise in computer algebra and managing three projects simultaneously, our learning curve was quite steep.

The primary tool for collaborative and version-controlled coding is \texttt{git}, so we needed to become familiar with sensible workflows. Due to initial indecision about whether or not to join the \texttt{OSCAR} system, we confusingly established different repositories with the same scope existing in parallel. If we had devised a plan ahead of coding, we would have been able to make informed decisions on what platform to use in which manner exactly, e.g., for website hosting.

It came as a great opportunity that we were able to join the \texttt{AlgebraicStatistics} module in \texttt{OSCAR}. As a consequence, we had a say in overarching design decisions, from the \texttt{OSCAR} codebase to the serialization system using the \texttt{.mrdi} file format. Additionally, we benefited from the existing technical and support structures. 
However, with all projects in active development, we had to adapt our internal progress to external progress.

We initially did not rely on a dedicated project-management tool (as we should have) but started multiple communication channels, effectively hindering efficient collaboration.
Typical hurdles such as fixed-term contracts and incompatible time zones were the straw that broke the camel's back, leading to a project runtime of nearly one and a half years, rather than the optimistically anticipated couple of months. In particular, MN continued to work for free for months after his student assistance officially ended (preparing himself exceptionally well for a career in academia).

\subsection{How to}
It is easy to criticize other projects, difficult to do better. {But} fellow mathematicians appreciate the effort! We thus wish for the development of \href{https://algebraicphylogenetics.org}{algebraicphylogenetics.org} to serve as a prototype for similar projects, leading to increasingly sustainable, usable, and trustworthy online mathematics. We recommend the following.

    \paragraph{Draw in all expertise available.} Discussions with data stewards, information specialists, computer scientists, computer algebraists, mathematicians who are in our community or own similar projects in other fields have helped us shape an idea of what is possible and what is needed to make algebraic-phylogenetics data both FAIR and reproducible.
    \paragraph{Consult research-data management experts.} In particular, we based our work on advice from research-data initiatives like \href{https://FAIRsFAIR.org}{FAIRsFAIR.org} and MaRDI as well as implementation discussions following the original publication of FAIR Principles. 
    \paragraph{Focus on usability.} All great work is for naught if it is not used. So we first compared notes with researchers who use Small Phylogenetic Trees actively, identifying pros and cons of the setup. We wondered who does not use the library, but would use it if we could make the content more understandable. As a result, we put the look-up table centrally on the new website's landing page to enable everyone to get a grasp of the content, with details for experts found further down the website tree.
    \paragraph{Investigate community needs.} The data layout plan and reviewing active research on algebraic phylogenetics enabled us to trace the workflow of researchers using Small Phylogenetic Trees. In order to help with the main hurdle, integrating results into their daily work, we put a priority on linking computational content and theory more directly.%
    \paragraph{Identify gaps.} The key points which make Small Phylogenetic Trees appear outdated are the early 2000s time stamps and partially broken list of maintainers. We want the Algebraic Phylogenetics Database to be up to date with current research, and visibly so, so we established contribution guidelines, versioning, and merged old data with new. To stop each researcher writing their own code to fill gaps between presented results and question of interest, we set up the collaborative software project.%
    \paragraph{Include documentation.} By making all code readable for humans, providing versioning of all digital resources, having permanent digital identifiers and interlinking as many objects as possible (models, papers, code etc.), we soften the impact of a number of dangers: that the software system we used may break, that the programming languages we use may die, that modern web technology evolves faster than we do; in short, that our currently very modern presentation will become one day as outdated as databases printed in books are today.
    
    \paragraph{Make it sustainable.} The job-interview question: where will our database be ten years from now? The math will still be true. To keep it FAIR and reproducible, we make it human as well as machine readable. We found an academic institution to host the website for the foreseeable future, and we established a community of maintainers with shared contact info, assigning one person the main responsibility. By choosing well-known \texttt{git} technology for external contributions, we enable the algebraic-phylogenetics community to actively sustain our project.
    \paragraph{Link it to larger infrastructure.} The FAIR Principles require metadata to be persistent, even if the data themselves are no longer available. To achieve this, we plan to put relevant metadata of the classification on Wikidata to supplement our Zenodo snapshots. This way, both data and metadata can be automatically harvested by online search engines and integrated into broader knowledge graphs for mathematics and neighboring disciplines. With the MaRDI Portal\footnote{\href{https://https://portal.mardi4nfdi.org}{portal.mardi4nfdi.org}} growing and interlinking formalized mathematics with paper theorems, the time for this type of integration is now.

\section*{Acknowledgements}

This project is a community effort.

We are hugely indebted to the authors of the original Small Phylogenetic Trees for enabling us to continue their legacy: Marta Casanellas, Serkan Ho\c{s}ten, Luis David Garcia Puente, Jacob Porter, Lior Pachter, Stacey Stokes, Bernd Sturmfels, and Seth Sullivant. In particular, we are very grateful to Luis David Garcia Puente, the principal maintainer, who provided us with valuable input in initial discussions, generous access to the original website files and data, and many thoughtful comments on this manuscript---and endured the pain we inflicted when criticizing his creation.

The software supplementing the database, \texttt{AlgebraicPhylogenetics.jl} and \texttt{AlgebraicStatistics.jl}, would not have been possible without the ongoing practical support and sharing of expertise of the \texttt{OSCAR} group: especially, Antony Della Vecchia, Lars Kastner, Michael Joswig, and Andrei Com\v{a}neci. We also wish to thank Tobias Boege and Ben Hollering for technical help and fruitful discussions about the implementation of mathematical concepts within our software project.

The Max Planck Institute for Mathematics in the Sciences supported MGL and MN in the early stages of this project and continues to support TB. We are immensely grateful for their acquisition of the domain and for hosting the \texttt{gitlab} fork of the GitHub project underlying \href{https://algebraicphylogenetics.org}{algebraicphylogenetics.org}.

TB, CG and MN were supported by the mathematical research-data initiative MaRDI, funded by the Deutsche Forschungsgemeinschaft (DFG), project number 460135501, NFDI 29/1 \enquote{MaRDI -- Mathematische Forschungsdateninitiative}. MGL was partially funded by the Institute for Computational and Experimental Mathematics (ICERM), under NSF Grant No. DMS-1929284, while in residence at the Fall 2024 Semester Program on Theory, Methods, and Applications of Quantitative Phylogenomics.

Pieter Belmans' \href{https://fanography.info}{fanography.info} served as a conceptual blueprint and inspiration for our project; we wish to thank him for making time to consult with us and for providing many helpful pointers before we embarked on this journey.







\bibliographystyle{apalike}
\bibliography{APD_refs.bib}
\end{document}